%% file: DFEffect.tex
\documentclass[12pt, leqno]{article}

\usepackage{latexsym,times,amsmath,amsfonts,amssymb,amsthm,bbm,natbib}

\input{DFEffect-Macros}

\begin{document}

\addtolength{\baselineskip}{+.44\baselineskip}
\leftmargini=0.4\parindent

\title{On Low-Dimensional Projections of High-Dimensional Distributions}
\author{Lutz D\"umbgen (University of Bern)\\
and\\
Perla Zerial (Technical University of Dresden)}
\date{July 2011, revised October 2011}
\maketitle

\begin{abstract}
Let $P$ be a probability distribution on $q$-dimensional space. The so-called Diaconis-Freed\-man effect means that for a fixed dimension $d << q$, most $d$-dimensional projections of $P$ look like a scale mixture of spherically symmetric Gaussian distributions. The present paper provides necessary and sufficient conditions for this phenomenon in a suitable asymptotic framework with increasing dimension $q$. It turns out that the conditions formulated by \cite{Diaconis_Freedman_1984} are not only sufficient but necessary as well. Moreover, letting $\hat{P}$ be the empirical distribution of $n$ independent random vectors with distribution $P$, we investigate the behavior of the empirical process $\sqrt{n}(\hat{P} - P)$ under random projections, conditional on $\hat{P}$.
\end{abstract}

\section{Introduction}

A standard method of exploring high-dimensional datasets is to examine various low-dimensional projections thereof. In fact, many statistical procedures are based explicitly or implicitly on a ``projection pursuit'', cf.\ \cite{Huber_1985}. \cite{Diaconis_Freedman_1984} showed that under weak regularity conditions on a distribution $P = P^{(q)}$ on $\R^q$, ``most'' $d$-dimensional orthonormal projections of $P$ are similar (in the weak topology) to a mixture of centered, spherically symmetric Gaussian distribution on $\R^d$ if $q$ tends to infinity while $d$ is fixed. A graphical demonstration of this disconcerting phenomenon is given by \cite{Buja_etal_1996}. Precise quantitative analyses are provided by Meckes (\nocite{Meckes_2009}{2009}, \nocite{Meckes_2011}{2011}) for situations where most projections are approximately Gaussian. The present paper provides further insight into the general phenomenon. We extend \nocite{Diaconis_Freedman_1984}{Diaconis and Freedman's~(1984)} results in two directions.

Section~\ref{sec:Population} gives necessary and sufficient conditions on the sequence $(P^{(q)})_{q \ge d}$ such that ``most'' $d$-dimensional projections of $P$ are similar to some distribution $Q$ on $\R^d$. It turns out that these conditions are essentially the conditions of \cite{Diaconis_Freedman_1984}. The novelty here is necessity. The limit distribution $Q$ is automatically a mixture of centered, spherically symmetric Gaussian distributions. The family of such measures arises in \cite{Eaton_1981} in a somewhat different context.

More precisely, let $\Gamma = \Gamma^{(q)}$ be uniformly distributed on the set of column-wise orthonormal matrices in $\R^{q \times d}$ (cf. Section~\ref{Proofs Population}). Defining
$$
	\gamma^\top P \ := \ \LL_{X \sim P}(\gamma^\top X)
$$
for $\gamma \in \R^{d \times q}$, we investigate under what conditions the random distribution $\Gamma^\top P$ converges weakly in probability to an arbitrary fixed distribution $Q$ as $q \to \infty$, while $d$ is fixed.

In Section~\ref{sec:Empirical} we study the relationship between $P = P^{(q)}$ and the empirical distribution $\hat{P} = \hat{P}^{(q,n)}$ of $n$ independent random vectors with distribution $P$, also independent from the projection matrix $\Gamma = \Gamma^{(q)}$. Suppose that the distributions $P^{(q)}$ satisfy the conditions of Section~\ref{sec:Population}. Then the random distributions $\hat{P}^{(q,n)}$ satisfy these conditions, too, as $q$ and $n$ tend to infinity. Furthermore, the standardized empirical measure $n^{1/2} \bigl( \Gamma^\top \hat{P} - \Gamma^\top P \bigr)$ satisfies a conditional Central Limit Theorem {\sl given the data} $\hat{P}$.

Proofs are deferred to Section~\ref{Proofs}. The main ingredients are \nocite{Poincare_1912}{Poincar\'{e}'s~(1912) Lemma} and a method invented by \cite{Hoeffding_1952} in order to prove weak convergence of conditional distributions. Further we utilize standard results from weak convergence and empirical process theory.

\section{The Diaconis-Freedman Effect}
\label{sec:Population}

Let us first settle some terminology. A random distribution $\hat{Q}$ on a separable metric space $(\M, \rho)$ is a mapping from some probability space into the set of Borel probability measures on $\M$ such that $\int f \, d \hat{Q}$ is measurable for any function $f \in \CC_b(\M)$, the space of bounded, continuous functions on $\M$. We say that a sequence $(\hat{Q}_k)_k$ of random distributions on $\M$ converges weakly in probability to some fixed distribution $Q$ if for each $f \in \CC_b(\M)$,
$$
	\int f \, d\hat{Q}_k \ \to_p \ \int f \, dQ \quad \mbox{as } k \to \infty .
$$
In symbols, $\hat{Q}_k \towp Q$ as $k \to \infty$. Standard approximation arguments (e.g.\ as in \nocite{vanderVaart_Wellner_1996}{van der Vaart and Wellner, 1996, Section~1.12}) show that $(\hat{Q}_k)_k$ converges in probability to $Q$ if, and only if,
$$
	D_{\rm BL}(\hat{Q}_k,Q) :=
	\sup_{f \in \FF_{\rm BL}} \, \Big| \int f \, d\hat{Q}_k - \int f \, dQ \Big| 
	\ \to_p \ 0 \quad (k \to \infty),
$$
where $\FF_{\rm BL}$ stands for the class of functions $f : \M \to [-1,1]$ such that $|f(x) - f(y)| \le \rho(x, y)$ for all $x,y \in \M$.

Now we can state the first result. Here and throughout, $\|\cdot\|$ denotes Euclidean norm and $\NN_{d,v}$ stands for the Gaussian distribution on $\R^d$ with mean vector $0$ and covariance matrix $v I_d$.

\begin{Theorem}
\label{thm:Diaconis-Freedman}
The following two assertions on the sequence $(P^{(q)})_{q \ge d}$ are equivalent:\\[1ex]
\textbf{(A1)} \ There exists a probability measure $Q$ on $\R^d$ such that
$$
	\Gamma^\top P \ \towp \ Q \quad\text{as} \ q \to \infty.
$$
\textbf{(A2)} \ If $X = X^{(q)}, \tilde{X} = \tilde{X}^{(q)}$ are independent random vectors with 
distribution $P$, then
$$
	\LL(\|X\|^2/q) \ \to_w \ R \quad \mbox{and} \quad X^\top\tilde{X}/q \ \to_p \ 0
	\quad \text{as} \ q \to \infty
$$
for some probability measure $R$ on $[0,\infty)$.

The limit distribution $Q$ in (A1) is a normal mixture, precisely,
$$
	Q \ = \ \int \NN_{d,v} \, R(dv)
$$
with the limiting distribution $R$ in (A2).
\end{Theorem}

\begin{Corollary}
\label{cor:Diaconis-Freedman}
The random probability measure $\Gamma^\top P$ converges weakly in probability to the standard Gaussian 
distribution $\NN_{d,1}$ if, and only if, the following condition is 
satisfied:\\[1ex]
\textbf{(B)} \ For independent random vectors $X = X^{(q)}, \tilde{X} = \tilde{X}^{(q)}$ with 
distribution $P$,
$$
	\|X\|^2/q \ \to_p \ 1 \quad \mbox{and} \quad X^\top\tilde{X}/q \ \to_p \ 0
	\quad\text{as} \ q \to \infty. \eqno{\Box}
$$
\end{Corollary}

The implication ``(A2) $\Longrightarrow$ (A1)'' in Theorem~\ref{thm:Diaconis-Freedman} as well as sufficiency of condition~(B) in Corollary~\ref{cor:Diaconis-Freedman} are due to \nocite{Diaconis_Freedman_1984}{Diaconis and Freedman~(1984, Theorem~1.1 and Proposition~4.2)}. They considered only (deterministic) empirical distributions $P$, but the extension to arbitrary distributions $P$ is straightforward; see also Section~\ref{sec:Empirical}.

It should be pointed out here that neither Theorem~\ref{thm:Diaconis-Freedman} nor Corollary~\ref{cor:Diaconis-Freedman} are just a consequence of \nocite{Poincare_1912}{Poincar\'{e}'s~(1912) Lemma}, although the latter is somehow at the heart of the proof. Poincar\'{e} showed that if $U_q = (U_{q,i})_{i=1}^q$ is uniformly distributed on the unit sphere in $\R^q$, then the Lebesgue density of $q^{1/2} U_{q,1}$ converges uniformly to the standard Gaussian density on $\R$. Translated into the present setting, one can show that for a fixed vector $x = x^{(q)} \in \R^q \setminus \{0\}$, the Lebesgue density of the random vector $\Gamma^\top x$ converges uniformly to the Lebesgue density of $\NN_{d,v}$ as $q \to \infty$ and $\|x\|^2/q \to v > 0$.

\begin{Example}
\label{Example 1}
Condition~(A2) is not a very restrictive requirement. For instance, suppose that $X = U (\mu_k + \sigma_k Z_k)_{k=1}^q$, where $(Z_k)_{k \ge 1}$ is a sequence of independent, identically distributed random variables with mean zero and variance one, while $U = U^{(q)}$ is an orthogonal matrix in $\R^{q\times q}$ and $\mu = \mu^{(q)} \in \R^q$, $\sigma = \sigma^{(q)} \in [0,\infty)^q$. Then condition~(A2) is satisfied if, and only if,
$$
	\|\mu\|^2/q \ \to \ 0, \quad 
	\|\sigma\|^2/q \ \to \ v \ge 0 \quad\text{and}\quad 
	\max_{1 \le k \le q} \, \sigma_k^2/q \ \to \ 0
	\leqno{\bf (A3)}
$$
as $q \to \infty$; see Section~\ref{Proofs}. Here $R = \delta_v$ and $Q = \NN_{d,v}$.
\end{Example}

\begin{Example}
\label{Example 2}
Suppose that $X \sim P^{(q)}$ has independent, identically distributed components such that
$$
	\Pr(X_i = \sqrt{q}) \ = \ 1 - \Pr(X_i = 0) \ = \ \pi_q ,
$$
where
$$
	\lim_{q \to \infty} \, q \pi_q \ = \ \lambda > 0 .
$$
Then $\LL(\|X\|^2/q) = \Bin(q, \pi_q) \to_w \mathrm{Poiss}(\lambda)$ and $\LL(X^\top \tilde{X}/q) = \Bin(q, \pi_q^2) \to_w \delta_0$ as $q \to \infty$. Hence (A2) is satisfied with $R = \mathrm{Poiss}(\lambda)$.
\end{Example}

\section{Empirical Distributions}
\label{sec:Empirical}

\paragraph{From $P$ to $\hat{P}$.}
If the distributions $P = P^{(q)}$ satisfy conditions~(A1-2), then the empirical distributions $\hat{P} = \hat{P}^{(q,n)}$ satisfy these conditions with high probability as $\min(q,n) \to \infty$. Precisely, one can easily deduce from condition~(A2) that
$$
	D_{\rm BL} \Bigl( \frac{1}{n} \sum_{i=1}^n \delta_{\|X_i\|^2/q}^{}, \, R \Bigr)
	\ \to_p^{} \ 0
$$
and
$$
	\frac{1}{n^2} \sum_{i,j=1}^n \min \bigl\{ |X_i^\top X_j^{}/q|, 1 \bigr\}
	\ \to_p^{} \ 0
$$
as $\min(q,n) \to \infty$. Thus Theorem~\ref{thm:Diaconis-Freedman} implies that
$$
	\Gamma^\top \hat{P}
	= \frac{1}{n} \sum_{i=1}^n \delta_{\Gamma^\top X_i}^{}
	\ \to_{w,p}^{} \ \int \NN_{d,v} \, R(dv)
$$
as both $q$ and $n$ tend to infinity, where the random projector $\Gamma$ and the empirical distribution $\hat{P}$ are assumed to be stochastically independent.

\paragraph{Comparing $P$ and $\hat{P}$, part~1.}
In some sense Theorem~\ref{thm:Diaconis-Freedman} is a negative, though mathematically elegant result. It warns us against hasty conclusions about high-dimensional data sets after examining a couple of low-dimensional projections. In particular, one should not believe in multivariate normality only because several projections of the data ``look normal''. On the other hand, even small differences between different low-dimens\-ional projections of $\hat{P}$ may be intriguing. Therefore we study the relationship between projections of the empirical distribution $\hat{P}$ and corresponding projections of $P$ in more detail.

In particular, we are interested in the halfspace norm
$$
	\|\Gamma^\top \hat{P} - \Gamma^\top P\|_{\rm KS} 
	\ := \ \sup_{\mbox{\small closed halfspaces } H \subset \R^d} \, 
		|\Gamma^\top \hat{P}(H) - \Gamma^\top P(H)|
$$
of $\Gamma^\top \hat{P} - \Gamma^\top P$. In case of $d = 1$ this is the usual Kolmogorov-Smirnov norm of $\Gamma^\top \hat{P} - \Gamma^\top P$. In what follows we use several well-known results from empirical process theory. Instead of citing original papers in various places we simply refer to the excellent monographs of \cite{Pollard_1984} and \cite{vanderVaart_Wellner_1996}. It is known that
\begin{equation}
	\Ex \sup_{\gamma \in \R^{q \times d}} \, \| \gamma^\top \hat{P} - \gamma^\top P\|_{\rm KS} 
	\ \le \ C \sqrt{q/n}
	\label{Alexander A}
\end{equation}
for some universal constant $C$. For the latter supremum is just the halfspace norm of $\hat{P} - P$, and generally the set of closed halfspaces in $\R^k$ is a Vapnik-Cervonenkis class with Vapnik-Cervonenkis index $k+1$. Inequality~(\ref{Alexander A}) does not capture the {\sl typical} deviation between $d$-dimensional projections of $\hat{P}$ and $P$. In fact,
$$
	\sup_{\gamma \in \R^{q \times d}} \, \Ex \| \gamma^\top \hat{P} - \gamma^\top P\|_{\rm KS} 
	\ \le \ C \sqrt{d/n} ,
$$
which implies that
\begin{equation}
	\Ex \| \Gamma^\top \hat{P} - \Gamma^\top P\|_{\rm KS} \ \le \ C \sqrt{d/n} .
	\label{Alexander B}
\end{equation}

Our next result implies the limiting distribution of $\sqrt{n} \| \Gamma^\top \hat{P} - \Gamma^\top P\|_{\rm KS}$ under conditions~(A1-2). More generally, let $\HH$ be a class of measurable functions from $\R^d$ into $[-1, 1]$. Any finite signed measure $M$ on $\R^d$ defines an element $h \mapsto M(h) := \int h \, dM$ of the space $\ell_\infty(\HH)$ of all bounded functions on $\HH$ equipped with supremum norm $\|z\|_\HH := \sup_{h \in \HH} |z(h)|$. We shall impose the following three conditions on the class $\HH$ and the distribution $Q = \int \NN_{d,v} \, R(dv)$:\\[1ex]
\textbf{(C1)} \ There exists a countable subset $\HH_o$ of $\HH$ auch that each $h \in \HH$ can be represented as pointwise limit of some sequence in $\HH_o$.\\[0.5ex]
\textbf{(C2)} \ The set $\HH$ satisfies the uniform entropy condition
$$
	\int_0^1 \sqrt{\log N(u, \HH)} \, du \ < \ \infty.
$$
Here $N(u, \HH)$ is the supremum of $N(u, \HH, \tilde{Q})$ over all probability measures $\tilde{Q}$ on $\R^d$, and $N(u, \HH, \tilde{Q})$ is the smallest number $m$ such that $\HH$ can be covered with $m$ balls having radius $u$ with respect to the pseudodistance
$$
	\rho_{\tilde{Q}}(g, h) \ := \ \sqrt{ \tilde{Q}((g - h)^2) }.
$$
\textbf{(C3)} \ For any sequence $(Q_k)_k$ of probability measures converging weakly to $Q$,
$$
	\|Q_k - Q\|_\HH \ \to \ 0 \quad \mbox{as } k \to \infty.
$$

Condition~(C1) ensures that random elements such as $\|\Gamma^\top \hat{P} - \Gamma^\top P\|_\HH$ are measurable. An example for conditions~(C1-2) is the set $\HH$ of (indicators of) closed halfspaces in $\R^d$. Then condition~(C3) is a consequence of general results by \cite{Billingsley_Topsoe_1967}, provided that $Q(\{0\}) = 0$, i.e.\ $R(\{0\}) = 0$.

A particular consequence of (C2) is existence of a centered Gaussian process $B_Q$, a so-called $Q$-bridge, having uniformly continuous sample paths with respect to $\rho_Q$ and covariances
$$
	\Ex \bigl( B_Q(g) B_Q(h) \bigr) \ = \ Q(gh) - Q(g) Q(h) ,
$$
which can be proved via a Chaining argument.

\begin{Theorem}
\label{thm:Empirical 1}
Suppose that the sequence $(P^{(q)})_{q \ge d}$ satisfies conditions~(A1-2) of Theorem~\ref{thm:Diaconis-Freedman}, and suppose that $\HH$ fulfills conditions~(C1-3). Then
$$
	B_{}^{(q,n)}
	\ := \ \Bigl( n^{1/2}
		\bigl( \Gamma^\top \hat{P} - \Gamma^\top P \bigr)(h)
	\Bigr)_{h \in \HH}
$$
converges in distribution in $\ell_\infty(\HH)$ to $B_Q$ as $\min(q,n) \to \infty$.
\end{Theorem}

\paragraph{Comparing $P$ and $\hat{P}$, part~2.}
Theorem~\ref{thm:Empirical 1} takes into account the randomness in both the data (i.e.\ $\hat{P}$) and the projection matrix $\Gamma$. However, exploratory projection pursuit means considering \textsl{several} projections of \textsl{one} data set. Thus we consider independent copies $\Gamma_\ell = \Gamma_\ell^{(q)}$, $\ell \ge 1$, of $\Gamma$ which are also independent from $\hat{P}$. With these projection matrices we define
$$
	B_\ell^{(q,n)}
	\ := \ \Bigl( n^{1/2}
		\bigl( \Gamma_\ell^\top \hat{P} - \Gamma_\ell^\top P \bigr)(h)
	\Bigr)_{h \in \HH}
$$
and study the distribution of
$$
	\bs{B}^{(q,n)} \ := \ \bigl( B_\ell^{(q,n)}(h) \bigr)_{(\ell,h) \in \Lambda \times \HH}
$$
for $\Lambda := \{1,\ldots,L\}$ with an arbitrary fixed integer $L \ge 1$.

Subsequently a particular decomposition of the $Q$-Brigde $B_Q$ will be used:
$$
	B_Q \ = \ B_Q' + B_Q''
$$
with stochastically independent and centered Gaussian processes $B_Q', B_Q''$ on $\HH$, where
\bea
	\Ex \bigl( B_Q'(g) B_Q'(h) \bigr)
	& = & Q(gh) - \int \NN_{d,v}(g) \, \NN_{d,v}(h) \, R(d v) \\
	& = & \int \bigl( \NN_{d,v}(gh) - \NN_{d,v}(g) \NN_{d,v}(h) \bigr)
		\, R(d v) \\
	\Ex \bigl( B_Q''(g) B_Q''(h) \bigr)
	& = & \int \NN_{d,v}(g) \NN_{d,v}(h) \, R(d v) - Q(g) Q(h) .
\eea
By means of \nocite{Anderson_1955}{Anderson's~(1955) Lemma} or a further application of Chaining one can show that both $B_Q'$ and $B_Q''$ admit versions with uniformly continuous sample paths.

\begin{Theorem}
\label{thm:Empirical 2}
Suppose that the conditions of Theorem~\ref{thm:Empirical 1} are satisfied. Further, let $B_{Q,1}', B_{Q,2}', B_{Q,3}', \ldots$ be independent copies of $B_Q'$ and independent from $B_Q''$. Then for any fixed integer $L \ge 1$, the process $\bs{B}^{(q,n)} = \bigl( B_\ell^{(q,n)}(h) \bigr)_{(\ell,h) \in \Lambda \times \HH}$ converges in distribution in $\ell_\infty(\Lambda \times \HH)$ to
$$
	\bs{B} \ := \ \bigl( B_{Q,\ell}'(h) + B_Q''(h) \bigr)_{(\ell,h) \in \Lambda \times \HH}
$$
as $\min(q,n) \to \infty$.
\end{Theorem}

\begin{Remark}[Understanding the decomposition $B_Q = B_Q' + B_Q''$ heuristically]
Note that $B^{(q,n)}(h) = \sqrt{n} \int h(\Gamma^\top x) \, (\hat{P} - P)(dx)$. Thus
\bea
	\Ex \bigl( B^{(q,n)}(h) \,\big|\, \hat{P} \bigr)
	& = & \sqrt{n} \int \Ex h(\Gamma^\top x) \, (\hat{P} - P)(dx) \\
	& = & \sqrt{n} \int \tilde{\NN}_{d,q,\|x\|}(h) \, (\hat{P} - P)(dx)
\eea
with $\tilde{\NN}_{d,q,\|x\|} := \LL(\Gamma^\top x)$. Here we utilize orthogonal invariance of $\LL(\Gamma)$. Consequently, $\Ex(B^{(q,n)} \,|\, \hat{P})$ is a standardized empirical process indexed by the special functions $x \mapsto \tilde{\NN}_{d,q,\|x\|}(h)$, $h \in \HH$, and
\bea
	\lefteqn{ \Ex \Bigl( \Ex \bigl( B^{(q,n)}(g) \,\big|\, \hat{P} \bigr)
		\Ex \bigl( B^{(q,n)}(h) \,\big|\, \hat{P} \bigr) \Bigr) } \\
	& = & \int \tilde{\NN}_{d,q,\|x\|}(g) \tilde{\NN}_{d,q,\|x\|}(h) \, P(dx)
		- \int \tilde{\NN}_{d,q,\|x\|}(g) \, P(dx)
			\int \tilde{\NN}_{d,q,\|x\|}(h) \, P(dx) .
\eea
Since $\tilde{\NN}_{d,q,\|x\|}$ is close to $\NN_{d,\|x\|^2/q}$ and $\LL(\|X\|^2/q)$ is close to $R$ for large $q$, the latter covariance is close to
$$
	\int \NN_{d,v}(g) \NN_{d,v}(h) \, R(dv)
		- \int \NN_{d,v}(g) \, R(dv)
			\int \NN_{d,v}(h) \, R(dv)
	\ = \ \Ex \bigl( B_Q''(g) B_Q''(h) \bigr) .
$$
\end{Remark}

\begin{Example}
Suppose that $d = 1$, and let $\HH$ consist of all indicator functions $1_{(-\infty,t]}$, $t \in \R$. Then Theorems~\ref{thm:Empirical 1} and \ref{thm:Empirical 2} are applicable whenever $R(\{0\}) = 0$. Writing $M(t)$ instead of $M(1_{(-\infty,t]})$, the covariance functions of $B_Q$, $B_Q'$ and $B_Q''$ are given by
\bea
	\Ex \bigl( B_Q(s) B_Q(t) \bigr)
	& = & Q(\min\{s,t\}) - Q(s) Q(t) , \\
	\Ex \bigl( B_Q'(s) B_Q'(t) \bigr)
	& = & Q(\min\{s,t\}) - \int \Phi(v^{-1/2} s) \Phi(v^{-1/2} t) \, R(dv) , \\
	\Ex \bigl( B_Q''(s) B_Q''(t) \bigr)
	& = & \int \Phi(v^{-1/2} s)\Phi(v^{-1/2} t) \, R(dv)
		- Q(s) Q(t)
\eea
for $s,t \in \R$, where $Q(u) = \int \Phi(v^{-1/2} u) \, R(dv)$, and $\Phi$ denotes the standard Gaussian distribution function.
\end{Example}

\begin{Remark}[Conservative inference]
Under conditions (A1-2) and (C1-3), \textsl{pretending} the empirical processes $B^{(q,n)}_\ell$, $1 \le \ell \le L$, to be independent and identically distributed leads typically to conservative procedures. Precisely, let $U$ be an open subset of $\ell_\infty(\HH)$. For instance let $U = \bigl\{ b \in \ell_\infty(\HH) : \|b\|_{\HH} < \kappa \bigr\}$ for some constant $\kappa > 0$. Then it follows from Theorem~\ref{thm:Empirical 2} that
$$
	\liminf_{\min(q,n) \to \infty}
		\, \Pr \bigl( B_\ell^{(q,n)} \in U \ \text{for} \ 1 \le \ell \le L \bigr)
	\ \ge \ \Pr(B_Q \in U)^L .
$$
This may be verified as follows: By Theorem~\ref{thm:Empirical 2} and the Portmanteau Theorem, the limes inferior on the left hand side is not smaller than
\bea
	\Pr \bigl( B_{Q,\ell}' + B_Q'' \in U \ \text{for} \ 1 \le \ell \le L \bigr)
	& = & \Ex \Pr \bigl( B_{Q,\ell}' + B_Q'' \in U \ \text{for} \ 1 \le \ell \le L
		\,\big|\, B_Q'' \bigr) \\
	& = & \Ex \Bigl(
		\Pr \bigl( B_Q' + B_Q'' \in U \,\big|\, B_Q'' \bigr)^L \Bigr) ,
\eea
and by Jensen's inequality the latter expression is not smaller than
$$
	\Bigl( \Ex \Pr \bigl( B_Q' + B_Q'' \in U \,\big|\, B_Q'' \bigr) \Bigr)^L \\
	\ = \ \Pr(B_Q' + B_Q'' \in U)^L
	\ = \ \Pr(B_Q \in U)^L .
$$

If (A.1-2) is strengthened to (B) and $\Pr(B_Q \in \partial U) = 0$, then the previous arguments lead to
$$
	\left.\begin{array}{c}
		\displaystyle
		\lim_{\min(q,n) \to \infty}
		\, \Pr \bigl( B_\ell^{(q,n)} \in U \ \text{for} \ 1 \le \ell \le L \bigr) \\
		\displaystyle
		\lim_{\min(q,n) \to \infty}
		\, \Pr \bigl( B_\ell^{(q,n)} \in \overline{U} \ \text{for} \ 1 \le \ell \le L \bigr)
	\end{array}\!\!\right\}
	\ = \ \Pr(B_Q \in U)^L ,
$$
because $B_Q'' \equiv 0$ almost surely.
\end{Remark}

\begin{Remark}[The conditional point of view]
Considering several projections of one data set means that we are interested in the \textsl{conditional} distribution of $n^{1/2}(\Gamma^\top \hat{P} - \Gamma^\top P)$, \textsl{given $\hat{P}$}. Indeed one may interpret Theorem~\ref{thm:Empirical 2} in the sense that for large $q$ and $n$,
$$
	\LL \bigl( B^{(q,n)} \,\big|\, \hat{P} \bigr)
	\ \approx \ \LL \bigl( B_Q' + B_Q'' \,\big|\, B_Q'' \bigr) .
$$
In case of the stronger condition~(B) in Corollary~\ref{cor:Diaconis-Freedman}, $B_Q'' \equiv 0$, and
$$
	\LL \bigl( B^{(q,n)} \,\big|\, \hat{P} \bigr)
	\ \approx \ \LL(B_Q) .
$$
Here are precise statements:

\begin{Corollary}
\label{cor:Empirical}
Suppose that the conditions of Theorem~\ref{thm:Empirical 1} are satisified. Let $F$ be any bounded and continuous functional on $\ell_\infty(\HH)$ such that $F(B^{(q,n)})$ is measurable for all $q \ge d$ and $n \ge 1$. Then
$$
	\Ex \bigl( F(B^{(q,n)}) \,\big|\, \hat{P} \bigr) 
	\ \to_{\LL} \ \Ex \bigl( F(B_Q' + B_Q'') \,\big|\, B_Q'' \bigr)
$$
as $\min(q,n) \to \infty$. In case of a degenerate distribution $R$,
$$
	\Ex \bigl( F(B^{(q,n)}) \,\big|\, \hat{P} \bigr) 
	\ \to_p \ \Ex F(B_Q)
$$
as $\min(q,n) \to \infty$.
\end{Corollary}
\end{Remark}

\section{Proofs}
\label{Proofs}

\subsection{Hoeffding's~(1952) trick}
\label{Hoeffding}

In connection with randomization tests, \cite{Hoeffding_1952} observed that weak convergence of 
conditional distributions of test statistics is equivalent to the weak convergence of the 
\textsl{unconditional} distribution of suitable statistics in $\R^2$. His result can be extended 
straightforwardly as follows.

\begin{Lemma}[Hoeffding]
\label{Hoeffding A}
For $k \ge 1$ let $X_k, \tilde{X}_k \in \X_k$ and $G_k \in \G_k$ be independent random variables, where $X_k, \tilde{X}_k$ are identically distributed. Further let $m_k$ be some measurable mapping from $\X_k \times \G_k$ into the separable metric space $(\M, \rho)$, and let $Q$ be a fixed Borel probability measure on $\M$. Then, as $k \to \infty$, the following two assertions are equivalent:
$$
	\LL \bigl( m_k(X_k, G_k) \,\big|\, G_k \bigr) \ \towp \ Q.
	\leqno{\bf(D1)}
$$
$$
	\LL \bigl( m_k(X_k, G_k), m_k(\tilde{X}_k, G_k) \bigr) \ \to_w \ Q \otimes Q.
	\leqno{\bf (D2)}
$$
\end{Lemma}

Applications of this equivalence with non-Euclidean spaces $\M$ are presented by \cite{Romano_1989}. We shall utilize Lemma~\ref{Hoeffding A} in order to prove Theorem~\ref{thm:Diaconis-Freedman}.

\begin{proof}[\bf Proof of Lemma~\ref{Hoeffding A}]
Define $Y_k := m_k(X_k, G_k)$ and $\tilde{Y}_k := m_k(\tilde{X}_k, G_k)$. Suppose first that (D2) ist true, i.e.\ $\LL(Y_k, \tilde{Y}_k) \to_w Q \otimes Q$. Then for any $f \in \CC_b(\M)$,
\bea
	\lefteqn{ \Ex \bigl( \bigl( \Ex(f(Y_k) \,|\, G_k) - Q(f) \bigr)^2 \bigr) } \\
	& = & \Ex \bigl( \Ex(f(Y_k) \,|\, G_k)^2 \bigr) - 2 Q(f) \, \Ex \Ex(f(Y_k) \,|\, G_k) 
		+ Q(f)^2 \\
	& = & \Ex \Ex \bigl( f(Y_k) f(\tilde{Y}_k) \,\big|\, G_k \bigr) 
		- 2 Q(f) \, \Ex \Ex(f(Y_k) \,|\, G_k) + Q(f)^2 \\
	& = & \Ex \bigl( f(Y_k) f(\tilde{Y}_k) \bigr) - 2 Q(f) \Ex f(Y_k) + Q(f)^2 \\
	& \to & \int f(y) f(\tilde{y}) \, Q(dy) Q(d\tilde{y}) - Q(f)^2 \\
	& = & 0 .
\eea
Thus $\LL(Y_k \,|\, G_k) \towp Q$.

On the other hand, suppose that (D1) is satisfied, i.e.\ $\LL(Y_k \,|\, G_k) \towp Q$. Then for arbitrary 
$f, g \in \CC_b(\M)$,
\bea
	\Ex \bigl( f(Y_k) g(\tilde{Y}_k) \bigr)
	& = & \Ex \Ex \bigl( f(Y_k) g(\tilde{Y}_k) \,\big|\, G_k \bigr) \\
	& = & \Ex \bigl( \Ex(f(Y_k) \,|\, G_k) \, \Ex(f(\tilde{Y}_k) \,|\, G_k) \bigr) \\
	& \to & Q(f) Q(g),
\eea
because $\Ex(h(Y_k) \,|\, G_k) \to_p \int h \, dQ$ and $\big| \Ex(h(Y_k) \,|\, G_k) \big| \le \|h\|_\infty < \infty$ for each $h \in \CC_b(\M)$. Thus we know that $\Ex F(Y_k, \tilde{Y}_k) \to \int F \, dQ \otimes Q$ for arbitrary functions $F(y, \tilde{y}) = f(y) g(\tilde{y})$ with $f, g \in \CC_b(\M)$. But this is known to be equivalent to weak convergence of $\LL(Y_k, \tilde{Y}_k)$ to $Q \otimes Q$; see \nocite{vanderVaart_Wellner_1996}{van der Vaart and Wellner~(1996, Chapter~1.4)}.

Here is an alternative argument: With $\hat{Q}_k := \LL(Y_k \,|\, G_k)$, Assumption (D1) is equivalent to $D_{\rm BL}(\hat{Q}_k,Q) \to_p 0$. To prove that $\LL(Y_k,\tilde{Y}_k) \to Q \otimes Q$, it suffices to show that $\Ex \bigl( F(Y_k,\tilde{Y}_k) \,\big|\, G_k \bigr) \to_p \int F \, d Q\otimes Q$ for any function $F : \M \times \M \to [-1,1]$ such that $\bigl| F(y,\tilde{y}) - F(z,\tilde{z}) \bigr| \le \rho(y,z) + \rho(\tilde{y},\tilde{z})$ for arbitrary $y, \tilde{y}, z, \tilde{z} \in \M$. But this entails that $F(y,\cdot), F(\cdot,\tilde{y}) \in \FF_{\rm BL}$ for arbitrary $y, \tilde{y} \in \M$. Consequently,
\bea
	\lefteqn{ \biggl| \Ex \bigl( F(Y_k,\tilde{Y}_k) \,\big|\, G_k \bigr) - \int F \, d Q \otimes Q
		\biggr| } \\
	& = & \biggl| \int F \, d \bigl( \hat{Q}_k \otimes \hat{Q}_k - Q \otimes Q \bigr) \biggr| \\
	& \le & \int \biggl| \int F(\cdot,\tilde{y}) \, d \bigl( \hat{Q}_k - Q \bigr) \biggr|
			\, \hat{Q}_k(d\tilde{y})
		+ \int \biggl| \int F(y,\cdot) \, d \bigl( \hat{Q}_k - Q \bigr) \biggr| \, Q(dy) \\
	& \le & 2 D_{\rm BL}(\hat{Q}_k,Q) .
\eea\\[-8ex]
\end{proof}

\subsection{Proofs for Section~\ref{sec:Population}}
\label{Proofs Population}

That $\Gamma = \Gamma^{(q)}$ is ``uniformly'' distributed on the set of column-wise orthonormal 
matrices in $\R^{q \times d}$ means that $\LL(U\Gamma) = \LL(\Gamma)$ for any fixed orthonormal 
matrix $U \in \R^{q \times q}$. For existence and uniqueness of the latter distribution we 
refer to \nocite{Eaton_1989}{Eaton~(1989, Chapters~1-2)}. For the present purposes the following explicit 
construction of $\Gamma$ described in \nocite{Eaton_1989}{Eaton~(1989, Chapter~7)} is sufficient. Let 
$Z = Z^{(q)} := (Z_1, Z_2, \ldots, Z_d)$ be a random matrix in $\R^{q \times d}$ with 
independent, standard Gaussian column vectors $Z_j \in \R^q$. Then
$$
	\Gamma \ := \ Z (Z^\top Z)^{-1/2}
$$
has the desired distribution, and
\begin{equation}
	\Gamma \ = \ q^{-1/2} Z \, (I + O_p(q^{-1/2})) \quad \mbox{as } q \to \infty.
	\label{Poincare}
\end{equation}
This equality can be viewed as an extension of \nocite{Poincare_1912}{Poincar\'{e}'s~(1912) Lemma}.

\begin{proof}[\bf Proof of Theorem~\ref{thm:Diaconis-Freedman}]
Let $\Gamma = \Gamma(Z)$ as above. Suppose that $Z = Z^{(q)}$, $X = X^{(q)}$ and $\tilde{X} = \tilde{X}^{(q)}$ are independent with $\LL(X) = \LL(\tilde{X}) = P$, and let $Y, \tilde{Y}$ be two independent random vectors in $\R^d$ with distribution $Q$. According to Lemma~\ref{Hoeffding A}, condition~(A1) is equivalent to
$$
	\begin{pmatrix} \Gamma^\top X \\ \Gamma^\top \tilde{X} \end{pmatrix} 
	\ \to_\LL \ \begin{pmatrix} Y \\ \tilde{Y} \end{pmatrix} .
	\leqno{\bf (A1')}
$$
Because of equation~(\ref{Poincare}) this can be rephrased as
$$
	\begin{pmatrix} Y^{(q)} \\ \tilde{Y}^{(q)} \end{pmatrix}
	:= \begin{pmatrix} q^{-1/2} Z^\top X \\ q^{-1/2} Z^\top \tilde{X} \end{pmatrix}
	\ \to_\LL \ \begin{pmatrix} Y \\ \tilde{Y} \end{pmatrix} .
	\leqno{\bf (A1'')}
$$
Now we prove equivalence of (A1'') and (A2) starting from the observation that
$$
	\LL \left( \begin{pmatrix} Y^{(q)} \\ \tilde{Y}^{(q)} \end{pmatrix} \right) 
	\ = \ \Ex \, \LL \left( \begin{pmatrix} Y^{(q)} \\ \tilde{Y}^{(q)} \end{pmatrix} 
		\,\Big|\, X, \tilde{X} \right) 
	\ = \ \Ex \, \NN_{2d}(0, \Sigma^{(q)}),
$$
where
$$
	\Sigma^{(q)} \ := \ \begin{pmatrix}
		q^{-1} \|X\|^2 \, I_d & q^{-1} X^\top \tilde{X} \, I_d \\
		q^{-1} X^\top \tilde{X} \, I_d & q^{-1} \|\tilde{X}\|^2 \, I_d
	\end{pmatrix} 
	\ \in \ \R^{2d \times 2d}.
$$

Suppose that condition~(A2) holds. Then $\Sigma^{(q)}$ converges in distribution to a random diagonal matrix
$$
	\Sigma \ := \ \begin{pmatrix} S^2 \, I_d & 0 \\ 0 & \tilde{S}^2 \, I_d \end{pmatrix}
$$
with independent random variables $S^2, \tilde{S}^2$ having distribution $R$. Clearly this implies that
$$
	\Ex \, \NN_{2d}(0, \Sigma^{(q)}) 
	\ \to_w \ \Ex \, \NN_{2d}(0, \Sigma) 
	\ = \ \LL \left( \begin{pmatrix} Y \\ \tilde{Y} \end{pmatrix} \right)
$$
with $Q = \Ex \, \NN_d(0, S^2 I_d)$. Hence (A1'') holds.

On the other hand, suppose that (A1'') holds. For any $t = (t_1^\top, t_2^\top)^\top \in \R^{2d}$, the Fourier transform of $\LL \bigl( (Y^{(q)}{}^\top, \tilde{Y}^{(q)}{}^\top)^\top \bigr)$ at $t$ equals
$$
	\Ex \, \exp \bigl( \bs{i} \, (t_1^\top Y^{(q)} + t_2^\top \tilde{Y}^{(q)}) \bigr) 
	\ = \ \Ex \, \exp(- t^\top \Sigma^{(q)} t/2) 
	\ = \ H^{(q)}(a(t)),
$$
where $\bs{i}$ stands for $\sqrt{-1}$, $a(t) := \bigl( \|t_1\|^2/2, \|t_2\|^2/2, t_1^\top t_2 \bigr)^\top \in \R^3$, and
$$
	H^{(q)}(a) 
	\ := \ \Ex \, \exp \bigl( - a_1 \|X\|^2/q - a_2 \|\tilde{X}\|^2/q 
		- a_3 X^\top \tilde{X}/q \bigr)
$$
denotes the Laplace transform of $\LL \bigl( \bigl( \|X\|^2/q, \|\tilde{X}\|^2/q, X^\top \tilde{X}/q \bigr)^\top \bigr)$ at $a \in \R^3$. By assumption, the Fourier transform at $t$ converges to
$$
	\Ex \exp(\bs{i} \, t_1^\top Y) \, \Ex \exp(\bs{i} \,  t_2^\top Y).
$$
Setting $t_2 = 0$ and varying $t_1$ shows that the Laplace transform of $\LL(\|X\|^2/q)$ converges pointwise on $[0,\infty)$ to a continuous function. Hence $\|X\|^2/q$ converges in distribution to some random variable $S^2 \ge 0$, and $Q = \Ex \NN_{d,S^2}$. Therefore, if $\tilde{S}^2$ denotes an independent copy of $S^2$, we know that $H^{(q)}(a(t))$ converges to
$$
	\Ex \, \exp(- a_1(t) S^2) \Ex \, \exp(- a_2(t) S^2) 
	\ = \ \Ex \, \exp \bigl( - a_1(t) S^2 - a_2(t) \tilde{S}^2 - a_3(t) \cdot 0 \bigr).
$$
A problem at this point is that for dimension $d = 1$ the set $\{a(t) : t \in \R^{2d}\} \subset \R^3$ has empty interior. Thus we cannot apply the standard argument about weak convergence and convergence of Laplace transforms. However, letting $t_2 = \pm t_1$ with $\|t_1\|^2/2 = 1$, one may conclude that
\bea
	0 & = & \lim_{q \to \infty} \bigl( 
		H^{(q)}(1, 1, 2) + H^{(q)}(1, 1, -2) - 2 H^{(q)}(1, 0, 0)^2 \bigr) \\
	& = & \lim_{q \to \infty} \bigl( 
		H^{(q)}(1, 1, 2) + H^{(q)}(1, 1, -2) 
		- 2 \, \Ex \exp(- \|X\|^2/q - \|\tilde{X}\|^2/q) \bigr) \\
	& = & 2 \, \lim_{q \to \infty} \,
		\Ex \Bigl( \exp \bigl( - \|X\|^2/q - \|\tilde{X}\|^2/q \bigr) 
			\bigl( \cosh(2 X^\top \tilde{X}/q) - 1 \bigr) \Bigr) .
\eea
But for arbitrary small $\eps > 0$ and large $r > 0$,
\bea
	\lefteqn{ \Ex \Bigl( \exp \bigl( - \|X\|^2/q - \|\tilde{X}\|^2/q \bigr) 
			\bigl( \cosh(2 X^\top \tilde{X}/q) - 1 \bigr) \Bigr) } \\
	& \ge & \exp(- 2 r) (\cosh(2 \eps) - 1) \,
		\Pr \bigl( \|X\|^2/q < r, \|\tilde{X}\|^2/q < r, |X^\top \tilde{X}/q| 
		\ge \eps \bigr) \\
	& \ge & \exp(- 2 r) (\cosh(2 \eps) - 1) \,
		\Bigl( \Pr \bigl( |X^\top \tilde{X}/q| \ge \eps \bigr)
			- 2 \Pr(\|X\|^2/q \ge r) \Bigr) \\
	& \ge & \exp(- 2 r) (\cosh(2 \eps) - 1) \,
		\Bigl( \Pr \bigl( |X^\top \tilde{X}/q| \ge \eps \bigr)
			- 2 \Pr(S^2 \ge r) + o(1) \Bigr) .
\eea
Hence
$$
	\limsup_{q \to \infty} \, \Pr \bigl( |X^\top \tilde{X}/q| \ge \eps \bigr) 
	\ \le \ 2 \Pr(S^2 \ge r) .
$$
Letting $r \to \infty$ shows that $X^\top \tilde{X}/q \to_p 0$.
\end{proof}

\begin{proof}[\bf Proof of equivalence of (A2) and (A3)]
Proving that (A3) implies (A2) is elementary. In order to show that (A2) implies (A3) note first that conditions~(A2) for the distributions $P^{(q)}$ imply the same conditions for the symmetrized distributions
$$
	P_o = P_o^{(q)} \ := \ \LL(X - \tilde{X}) 
	\ = \ \LL \Bigl( \bigl( \sigma_k(Z_k - Z_{q+k}) \bigr)_{1 \le k \le q} \Bigr).
$$
Condition~(A2) for these distributions reads as follows.
\begin{eqnarray}
	\LL \Bigl( \sum_{k = 1}^q (Z_k - Z_{q+k})^2 \sigma_k^2/q \Bigr) 
	& \to_w & R_o = R \star R \quad\text{and} 
	\label{A3a} \\
	\sum_{k = 1}^q (Z_k - Z_{q+k}) (Z_{2q+k} - Z_{3q+k}) \sigma_k^2/q 
	& \to_p & 0.
	\label{A3b}	
\end{eqnarray}
The factors $(Z_k - Z_{q+k}) (Z_{2q+k} - Z_{3q+k})$, $1 \le k \le q$, in 
(\ref{A3b}) are independent, identically and symmetrically distributed. By conditioning on any one of these factors one can deduce from (\ref{A3b}) that $\max_{1 \le k \le q} \sigma_k^2/q \to 0$. But then
$$
	\sum_{k = 1}^q \sigma_k^2 (Z_k - Z_{q+k})^2/q 
	\ = \ 2 \|\sigma\|^2/q + o_p(1 + \|\sigma\|^2/q),
$$
and one can deduce from (\ref{A3a}) that $\|\sigma\|^2/q$ converges to some fixed 
number $v$; in particular, $R = \delta_v$. Now we return to the original distributions $P$. Here 
the second half of (A2) means that
\bea
	\lefteqn{ \sum_{k = 1}^k (\mu_k + \sigma_k Z_k)(\mu_k + \sigma_k Z_{q + k})/q } \\
	& = & \|\mu\|^2/q + \sum_{k = 1}^q \mu_k \sigma_k (Z_k + Z_{q+k}) /q
		+ \sum_{k = 1}^q \sigma_k^2 Z_k Z_{q+k}/q \\
	& = & o_p(1).
\eea
Since
\bea
	\Ex \biggl( \Bigl( \sum_{k = 1}^q \mu_k \sigma_k (Z_k + Z_{q+k})/q \Bigr)^2 \biggr) 
	& = & \sum_{k = 1}^q \mu_k^2 \sigma_k^2 / q^2 \ = \ o(\|\mu\|^2/q), \\
	\Ex \biggl( \Bigl( \sum_{k = 1}^q \sigma_k^2 Z_k Z_{q+k}/q \Bigr)^2 \biggr) 
	& = & \sum_{k = 1}^q \sigma_k^4/q^2 \ \to \ 0,
\eea
it follows that $\|\mu\|^2/q \to 0$.
\end{proof}

\subsection{Proofs for Section~\ref{sec:Empirical}}

Since Theorem~\ref{thm:Empirical 1} is just Theorem~\ref{thm:Empirical 2} with $L = 1$, it suffices to verify the latter.

\begin{proof}[\bf Proof of Theorem~\ref{thm:Empirical 2}]
It suffices to verify the following two 
claims:\\[1ex]
\textbf{(F1)} \ As $q \to \infty$ and $n \to \infty$, the finite-dimensional marginal distributions of the process $\bs{B}^{(q,n)}$ converge to the corresponding finite-dimensional distributions of $\bs{B}$.\\[0.5ex]
\textbf{(F2)} \ As $q \to \infty$, $n \to \infty$ and $\delta \downarrow 0$,
$$
	\max_{\ell \in \Lambda} \ \sup_{g, h \in \HH : \rho_Q(g, h) < \delta} \, 
		\Big| B^{(q,n)}_\ell(g) - B^{(q,n)}_\ell(h) \Big| \ \to_p \ 0.
$$

The second condition, (F2), means that the processes $\bs{B}^{(q,n)}$ are asymptotically equi\-continuous with respect to the pseudodistance
$$
	\bs{\rho}_Q \bigl( (\ell,g), (m,h) \bigr)
	\ := \ 1\{\ell \ne m\} + \rho_Q(g,h)
$$
on $\Lambda \times \HH$.

In order to verify assertions~(F1-2) we consider the conditional distribution of $\bs{B}^{(q,n)}$ given the random matrix
$$
	\bs{\Gamma} = \bs{\Gamma}^{(q)} \ := \ (\Gamma_1, \Gamma_2, \ldots, \Gamma_L) \ \in \ \R^{q \times Ld}.
$$
In fact, if we define
$$
	f_{\ell,h}(\bs{v}) \ := \ h(v_\ell)
	\quad\text{for} \ 
	\bs{v} = (v_1^\top, \ldots, v_L^\top)^\top \in \R^{Ld},
$$
then
$$
	B^{(q,n)}_\ell(h) \ = \ n^{1/2}(\bs{\Gamma}^\top \hat{P} - \bs{\Gamma}^\top P)(f_{\ell,h}).
$$
Thus $\LL(\bs{B}^{(q,n)} \,|\, \bs{\Gamma})$ is essentially the distribution of an empirical process based on $n$ independent random vectors with distribution $\bs{\Gamma}^\top P$ on $\R^{Ld}$ and indexed by the family $\tilde{\HH} := \{f_{\ell,h} : \ell \in \Lambda, h \in \HH\}$.

The multivariate version of Lindeberg's Central Limit Theorem entails that for large $q$ and $n$, the finite-dimensional marginal distributions of $\bs{B}^{(q,n)}$, conditional on $\bs{\Gamma}$, can be approximated by the corresponding finite-dimensional distributions of a centered Gaussian process on $\Lambda \times \HH$ with the same covariance function, namely,
\bea
	\Sigma^{(q)} \bigl( (\ell, g), (m, h) \bigr) 
	& := & \Cov \bigl( B^{(q,n)}_\ell(g), B^{(q,n)}_m(h) \,\big|\, \bs{\Gamma} \bigr) \\
	& = & \bs{\Gamma}^\top P(f_{\ell,g} f_{m,h}) 
		- \bs{\Gamma}^\top P(f_{\ell,g}) \bs{\Gamma}^\top P(f_{m,h}).
\eea
It follows from equality~(\ref{Poincare}) and the proof of Theorem~\ref{thm:Diaconis-Freedman} that
$$
	\bs{\Gamma}^\top P \ \towp \ \bs{Q} \ := \ \int \NN_{Ld,v} \, R(dv) 
	\quad \text{as} \ q \to \infty ,
$$
and this should imply convergence of $\Sigma^{(q)}$ to some limiting function as well. It was shown by \cite{Billingsley_Topsoe_1967} that condition~(C3) is equivalent to
\begin{equation}
	\lim_{\delta \downarrow 0} \ \sup_{h \in \HH} \,
	Q \Bigl\{ y \in \R^d : 
		\sup_{z : \|z - y\| < \delta} \, |h(z) - h(y)| > \eps \Bigr\} \ = \ 0 
	\quad \mbox{for any } \eps > 0.
	\label{C3'}
\end{equation}
Note that the $d$-dimensional marginal distributions of $\bs{Q}$ are just $Q$. Therefore one 
can easily deduce from (\ref{C3'}) that for any fixed $\eps > 0$,
$$
	\lim_{\delta \downarrow 0} \ \sup_{f', f'' \in \tilde{\HH} \cup \{1\}} \, 
		\bs{Q} \Bigl\{ \bs{v} \in \R^{Ld} : \sup_{\bs{w} : \|\bs{w} - \bs{v}\| < \delta} \, 
			|f'f''(\bs{w}) - f'f''(\bs{v})| > \eps \Bigr\} 
	\ = \ 0 .
$$
Hence a second application of \cite{Billingsley_Topsoe_1967} shows that
\begin{equation}
	\sup_{f',f'' \in \tilde{\HH} \cup \{1\}} \, 
		|\bs{\Gamma}^\top P (f'f'') - \bs{Q}(f'f'')| \ \to \ 0 
	\quad\text{as} \ q \to \infty ,
	\label{Topsoe}
\end{equation}
because $\bs{\Gamma}^\top P \to_{w,p} \bs{Q}$. In particular, the conditional covariance function $\Sigma^{(q)}$ converges uniformly in probability to the covariance function $\Sigma$, where
\bea
	\lefteqn{ \Sigma \bigl( (\ell, g), (m, h) \bigr) 
		\ := \ \bs{Q}(f_{\ell,g} f_{m,h}) 
		- \bs{Q}(f_{\ell,g}) \bs{Q}(f_{m,h}) } \\
	& = & \int \NN_{Ld,v}(f_{\ell,g} f_{m,h}) \, R(dv) - Q(g) Q(h) \\
	& = & \begin{cases}
		\displaystyle
		\int \NN_{d,v}(gh) \, R(dv)  - Q(g) Q(h)
			& \text{if} \ \ell = m, \\[1.5ex]
		\displaystyle
		\int \NN_{d,v}(g) \NN_{d,v}(h) \, R(dv)  - Q(g) Q(h) 
			& \text{if} \ \ell \neq m,
		\end{cases} \\[0.5ex]
	& = & \Cov \bigl( B_{Q,\ell}'(g) + B_Q''(g), B_{Q,m}'(h) + B_Q''(h) \bigr)
\eea
as $q \to \infty$. This proves assertion~(F1).

As for assertion~(F2), it is well-known from empirical process theory that conditions~(C1-2) 
imply that for arbitrary fixed $\eps > 0$,
\begin{equation}
	\max_{\ell \in \Lambda} \,
	\Pr \Bigl( \sup_{g, h \in \HH : \rho_\ell^{(q)}(g, h) < \delta} \, 
		\Big| B_\ell^{(q,n)}(g) - B_\ell^{(q,n)}(h) \Big| \ge \eps \,\Big|\, \bs{\Gamma} \Bigr) 
	\ \to_p \ 0
	\label{almost F2}
\end{equation}
as $\min(q,n) \to \infty$ and $\delta \downarrow 0$. Here
$$
	\rho_\ell^{(q)}(g, h) 
	\ := \ \sqrt{ \bs{\Gamma}^\top P((f_{\ell,g} - f_{\ell,h})^2) }
	\ = \ \sqrt{ \Gamma_\ell^\top P((g - h)^2) } .
$$
But it follows from (\ref{Topsoe}) that
$$
	\max_{\ell \in \Lambda} \,
	\sup_{g, h \in \HH} \,
		|\rho_\ell^{(q)}(g, h)^2 - \rho_Q(g, h)^2| \ \to_p \ 0 
$$
as $q \to \infty$. Hence one may replace $\rho_\ell^{(q)}$ in (\ref{almost F2}) with $\rho_Q$ and obtain assertion~(F2).
\end{proof}

\begin{proof}[\bf Proof of Corollary~\ref{cor:Empirical}]
The main trick is to replace conditional expectations with suitable sample means. Note that conditional on $\hat{P}$, the processes $B_1^{(q,n)}, B_2^{(q,n)}, B_3^{(q,n)}, \ldots$ are independent copies of $B_{}^{(q,n)}$. Likewise, conditional on $B_Q''$, the processes $B_{Q,1}' + B_Q'', B_{Q,2}' + B_Q'', B_{Q,3}' + B_Q'', \ldots$ are independent copies of $B_Q' + B_Q''$. Hence
$$
	\left. \begin{array}{c}
		\displaystyle
		\Ex \, \Bigl| \Ex \bigl( F(B^{(q,n)}) \,\big|\, \hat{P} \bigr)
		- L^{-1} \sum_{\ell=1}^L F(B^{(q,n)}_\ell) \Bigr| \\
		\displaystyle
		\Ex \, \Bigl| \Ex \bigl( F(B_Q' + B_Q'') \,\big|\, B_Q'' \bigr)
		- L^{-1} \sum_{\ell=1}^L F(B_{Q,\ell}' + B_Q'') \Bigr|
	\end{array} \!\! \right\}
	\ \le \ L^{-1/2} \|F\|_\infty
$$
for any integer $L \ge 1$. Consequently it suffices to show that for any fixed $L \ge 1$, the random variable $L^{-1} \sum_{\ell=1}^L F(B^{(q,n)}_\ell)$ converges in distribution to the random variable $L^{-1} \sum_{\ell=1}^L F(B_{Q,\ell}' + B_Q'')$ as $\min(q,n) \to \infty$. But this is a consequence of Theorem~\ref{thm:Empirical 2} and the Continuous Mapping Theorem, because
$$
	\bs{b} = \bigl( b_\ell(h) \bigr)_{(\ell,h) \in \Lambda\times\HH}
	\ \mapsto \ L^{-1} \sum_{\ell=1}^L F(b_\ell)
$$
defines a continuous mapping from $\ell_\infty(\Lambda\times\HH)$ to $\R$.
\end{proof}

\paragraph{Acknowledgement.}
Part of this work is contained in the diploma thesis of Perla Zerial (1995, Univ.\ of Heidelberg). We are grateful to Jon Wellner, Aad van der Vaart and an anonymous referee for their interest in this work, stimulating discussions and pertinent questions.

\end{document}

%% file: DFEffect-Macros.tex
\newtheoremstyle{localthm}
	{5pt} 
	{5pt} 
	{\sl} 
	{} 
	{\bf} 
	{} 
	{.7em} 
	{} 

\newtheoremstyle{localrem}
	{5pt} 
	{5pt} 
	{\rm} 
	{} 
	{\bf} 
	{} 
	{.7em} 
	{} 

\theoremstyle{localthm}
\newtheorem{Theorem}{Theorem}[section]
\newtheorem{Lemma}[Theorem]{Lemma}
\newtheorem{Corollary}[Theorem]{Corollary}

\theoremstyle{localrem}
\newtheorem{Remark}[Theorem]{Remark}
\newtheorem{Example}[Theorem]{Example}

\def\bea{\begin{eqnarray*}}
\def\eea{\end{eqnarray*}}

\def\towp{\to_{w,p}}

\def\bs{\boldsymbol}

\def\hat{\widehat}

\def\G{\mathbb{G}}
\def\M{\mathbb{M}}
\def\R{\mathbb{R}}
\def\X{\mathbb{X}}

\def\CC{\mathcal{C}}
\def\FF{\mathcal{F}}
\def\HH{\mathcal{H}}

\def\Ex{\mathop{\mathrm{I\!E}}\nolimits}
\def\Pr{\mathop{\mathrm{I\!P}}\nolimits}

\def\Bin{\mathrm{Bin}}
\def\Cov{\mathop{\mathrm{Cov}}\nolimits}

\def\LL{\mathcal{L}}
\def\NN{\mathcal{N}}

\def\eps{\epsilon}